\documentclass[leqno, 11pt]{article}
\usepackage{color}
\usepackage{graphicx}

\begin{document}

\begin{center}
{ \bf 
ON THE GROWTH RATES FOR  

A THREE-LAYER
  FLOW IN POROUS MEDIA }
 
 \hspace{0.25cm}
 
 GELU I. PA\c{S}A
 

\end{center}

\noindent
{\small 
We  study the  
displacement of three
 immiscible Stokes fluids
with constant viscosities
in a porous medium.
The middle-layer fluid is
contained in a bounded region.
We give an analysis of the linear stability of this process.
 This stability problem has no solution (in general).
 }

 \vspace{0.25cm}



\begin{center}
1. INTRODUCTION
\end{center}

The flow porous media are 
often studied by using 
the Hele-Shaw model. This is  obtained
by an average procedure of 
 a Stokes flow
in the narrow gap 
between two parallel plates.
The ﬂuid velocity is related to pressure gradients
via a Darcy-type law  
- see \cite{BE},
\cite{HS}.

The Saffman-Taylor instability  
\cite{SAFF-TAY} appears
when 
a less viscous 
Stokes fluid is displacing a more viscous one in a two-layer Hele-Shaw cell.
It is  important to find flow-models which could minimize this instability.  
A possible strategy is to
put some intermediate
liquids between the initial fluids. 
In \cite{PA}, \cite{PA2} we proved that  this strategy is not giving  us  an almost stable flow, even if the 
number of intermediate layers 
is very large.

In this paper  
we study  the linear stability
of the flow with 
a single intermediate liquid layer with a constant 
viscosity, first studied in
 \cite{DA0}.
 The growth constants of 
perturbations are the
eigenvalues of the stability system and are contained in 
the boundary conditions,
then we have also a compatibility condition.
 The eigenfunctions are the amplitudes
 of the linear perturbations, and   must be bounded (only small perturbations are allowed). 
 From the compatibility 
 condition it results
  some (unexpected) restrictions on the three viscosities,
which can lead to physical contradictions. 
 Thus the linear stability problem has no solution (in general).

{\begin{center}

  2. THE 3-LAYER 
  FLOW
  WITH CONSTANT  
  VISCOSITIES
  
\end{center}

We consider a horizontal 
Hele-Shaw cell in the fix plane $x_1Oy$, filled by three immiscible Stokes fluids
with constant viscosities.
The averaged
velocity $(u,v)$ 
is related with the pressure
 $p$ by the Darcy-type law 
 (\ref{A2}).

The displacing fluid is pushing
 the middle-layer fluid
with  the far-upstream 
velocity  $(U,0$).   
 The middle-layer fluid
 contained in a bounded region
is pushing the displaced liquid (say, oil).
This model is described in paper  \cite{DA0}. 
We use  the moving reference $x=x_1-Ut$.
The flow is governed
by the equations 
(1a), (1b) of 
\cite{DA0}:
\begin{equation}\label{A1}
 \nabla \cdot {\bf u}=0, 
 \quad {\bf u}=(u,v),
\end{equation}
\begin{equation}\label{A2}
 \nabla p = - \nu {\bf u}.
\end{equation}
Here the  viscosity $\nu$
verifies the conditions
$$
\nu=\mu, \, x \in (a, b); \quad \quad 
0< \mu_L < \mu < \mu_R;  \quad \quad b \leq 0;                      $$
\begin{equation}\label{A3}
\nu = \mu_L, \, x<a; \quad
\nu = \mu_R, \, x>b.
\end{equation}
The system 
(\ref{A1})-(\ref{A3})
is stationary along with
two planar interfaces 
$x=a, \,\, x=b$
separating these three fluid layers -
see \cite{GOR-HOM-1}.
On $x=a$, $x=b$   we 
assume the existence of the  positive  surface tensions 
$T(a),T(b)$. The Laplace-Young law will be  used on both interfaces.

\begin{center}

 3. THE LINEAR STABILITY ANALYSIS

\end{center}

 The perturbations  of 
 velocity and pressure
are denoted
 by $u', v', p'$.
 We consider
 a small positive number $\varepsilon$, 
 the wavenumbers $k$, the growth rates $\sigma$,
 the amplitude $f$ and
 the Fourier decomposition
$$  u'= \varepsilon
f(x)\exp(iky + \sigma t),
                      $$
\begin{equation}\label{J2}
f(x) = f(a) e^{  k(x-a) }, \,\, \forall x \leq a; \quad 
f(x) = f(b) e^{ -k(x-b) }, \,\, \forall x \geq b.
\end{equation}
The last condition 
is imposed   in 
\cite{SAFF-TAY}:
 $f$  must decay to zero in the far field and is continuous -
see also \cite{GOR-HOM-1}.   
The amplitude $f$ must 
be bounded also 
in the intermediate region.
Only small perturbations
are allowed in the frame of the linear perturbations. Thus 
we assume 
\begin{equation}
\label{J2A}
- \infty < 
f(k,x) < \infty, \quad 
\forall   x \in [a,b],
\quad \forall k \geq 0.
\end{equation}
The restriction 
(\ref{J2A})
is not considered in  \cite{DA0}.

We insert the perturbations
in equations 
(\ref{A1})-(\ref{A2}).
The 
linearized disturbance 
equations
and also   the linearized dynamic and kinematic interfacial conditions given  in  \cite{DA0}, \cite{DP} are used. 
Consider
$a=-L < b \leq 0$.
We get the problem (2a), (2b), (3)  of  \cite{DA0}
 (see also the relation (20) in \cite{DA1}): 
 
 \newpage

 \begin{equation}\label{A5}
 f_{xx}-k^2 f=0, \quad 
 x  \in (a,b), \quad 
 \forall k \geq 0;
\end{equation}
\begin{equation}\label{A6}
 f_x(a)= 
 [-E_a(k)/ \sigma +s] f(a);
\end{equation}
\begin{equation}\label{A7}
 f_x(b)= 
[E_b(k) / \sigma  +q]f(b);
\end{equation}
$$
 E_b(k):=\frac{
 (\mu_R-\mu)Uk^2-
 T(b)k^4}{\mu},  \quad
 q= -\frac{\mu_Rk}{\mu};
                       $$
\begin{equation}\label{A9}
 E_a(k): =  \frac{
 (\mu-\mu_L)Uk^2
 -  T(a)k^4}{\mu} , \quad
 s= \frac{\mu_Lk}{\mu}.
\end{equation}
 The same growth constants 
 appear in the boundary conditions, thus:
\begin{equation}\label{COMP}
\sigma(k)= 
\frac{kE_b(k)f(b)}
{\mu f_x^-(b)+
\mu_R kf(b)}      
= \frac{kE_a(k)f(a)}
{\mu_Lkf(a)- 
\mu f_x^+(a)},   \quad\forall k  \geq 0.
\end{equation}
Here $f_x^-(b), f_x^+(a)$ are the limit values of
$f_x$ in the points $b,a$.
{\it This condition is not mentioned  in   \cite{DA0}. }

\vspace{0.25cm}

The possible  
solutions of (\ref{A5})
(considered also in 
\cite{DA0}) are 
\begin{equation}\label{D1}
 f(x)= A(k) e^{kx} + B(k)e^{-kx}.
\end{equation}
In fact, in the cited paper 
is only specified that 
the 
$A, B$ are constants. 
We will often use the notation $A,B$ 
instead of  $A(k),B(k)$. {\it This is a new element of our paper}.
The condition
(\ref{COMP}) is equivalent with  the relations 
(\ref{F-4}), (\ref{F-5})
below
\begin{equation}\label{F-4}
E_b/E_a =
[\mu F(k,b)+ \mu_R] /
[\mu_L - \mu F(k,a)],
\end{equation}
\begin{equation}\label{F-5}
 F(k,x)= 
 \frac{Ae^{kx}-Be^{-kx}}
      {Ae^{kx}+Be^{-kx}}.
\end{equation}

From the boundary conditions
(\ref{A6})-(\ref{A7}) we obtain the relations
$$ 
[\mu f_x(a)-\mu_L f(a)] = -f(a) E_a/\sigma, \quad
[\mu f_x(b)+ \mu_R f(b)]
= f(b) E_b/\sigma .             $$
Therefore  $A,B$, 
are verifying the 
equations 
$$ \mu (Ae^{ka}-Be^{-ka})=
(-E_a/\sigma + \mu_L)
(Ae^{ka}+Be^{-ka}),
                    $$
$$ \mu (Ae^{kb}-Be^{-kb})=
(E_b/\sigma - \mu_R)
(Ae^{kb}+Be^{-kb}),
                    $$ 
 $$ 
 Ae^{ka}(\mu-\mu_L 
 +E_a/\sigma) + 
 Be^{-ka}(- \mu-\mu_L
 +E_a/\sigma) =0,   $$
$$ Ae^{kb}(\mu+\mu_R
-E_b/\sigma)  +
Be^{-kb}(-\mu+\mu_R
-E_b/\sigma) =0,        $$ 
and we get the homogeneous system 
$$ 
 Ae^{ka} c + Be^{-ka} d=0,
 \quad        
 Ae^{kb} g + Be^{-kb} h=0;
                         $$
 $$ 
 c=(\mu_L-\mu-E_a/\sigma),\quad 
d=(\mu_L+\mu-E_a/\sigma),
                          $$ 
  \begin{equation}\label{F3A}
g=(\mu_R+\mu-E_b/\sigma),\quad 
h=(\mu_R-\mu-E_b/\sigma).
\end{equation}
 
A solution 
$(A,B) \neq (0,0)$ exist if the following conditions is verified 
\begin{equation}\label{F-3}
 e^{2k(a-b)}ch- gd=0,
 \quad \forall k \geq 0.
\end{equation}

It is difficult to verify whether the compatibility
condition 
(\ref{F-4}) is fulfilled. 
To this end, 
we derive some properties 
of the  growth rates 
for   large $ k $.
The numerical values of 
$\sigma(k)$
are obtained in \cite{DA0}, but only for
$ k \leq 3.5 $. 

There are an infinity of solutions for the system 
(\ref{F3A}). For an arbitrary $B(k)$ we get  $A(k)$ from   (\ref{F3A}), due to the relationship (\ref{F-3}).

\vspace{0.5cm} 

{\bf Proposition  1.}
{\it 
  From   (\ref{F-3}) 
 we get
$$ \lim_{k\to\infty}
\frac{E_b(k)}{ \sigma(k)} 
= \mu_R + \mu, \quad 
\lim_{k\to\infty}
\frac{E_b(k)}{ \sigma(k)} 
= \mu_L + \mu.
                       $$
Proof. We
introduce the notations 
$Q, \Delta, i,j,m, n$
below and (\ref{F-3}) gives us 
$$
\sigma^2(k) (Q ij - mn)+
\sigma(k)[ (mE_b+nE_a)
-Q(iE_b+jE_a)]          
                      $$
\begin{equation}\label{R1}
+
(QE_aE_b-E_aE_b)=0; \quad
Q= e^{2k(a-b)}.
\end{equation}
$$
i=\mu_L-\mu; \,\,j=\mu_R-\mu;\,\,
m=\mu_L+\mu; \,\,n=\mu_R+
\mu;                   $$
$$
\,\, \Delta :=(mE_b+nE_a)^2 -
2Q(mE_b+nE_a)(iE_b+jE_a)
                       $$
\begin{equation}\label{R2}
+ Q^2(iE_b+jE_a)^2 - 4
(Q ij - mn)
(QE_aE_b-E_aE_b).
\end{equation}
For large values of $k$ we have 
$$ 
e^{-2k}k^3  \approx 0,
\,\, e^{-2k} << k^3;  
\quad
e^{-2k}k^6  \approx 0,
\,\, e^{-2k}k^6  << k^6 
\Rightarrow            $$
$$ (mE_b+nE_a)
-Q(iE_b+jE_a) \approx (mE_b+nE _a),          $$
$$
(Q ij-mn) \approx -mn, \quad 
 QE_aE_b-E_aE_b \approx 
-E_aE_b,               $$
$$ 2Q(mE_b+nE_a)(iE_b+jE_a) \approx  0,           \quad 
Q^2(iE_b+jE_a)^2 
\approx  0,             $$ 
$$
\Delta \approx  
(mE_b+nE_a)^2
-4mn E_aE_b=
(mE_b-nE_a)^2, 
                       $$ 
\begin{equation}\label{R3}
\sigma(k) \approx 
\frac{-(mE_b+nE_a) 
\,\,\, ^+_- \,\,\,
  |mE_b-nE_a|}{-2mn}.   
\end{equation}
We see that
$E_a(k), E_b(k) \neq 0$
if $k$ is large enough.
So  we have:
\begin{equation}\label{R4}
\lim_{k\to\infty}
\frac{E_b(k)}{ \sigma(k)} 
=  {\mu_R+\mu}, \,\,
 \lim_{k\to\infty} g(k) = 0, \,\,
 \lim_{k\to\infty} h(k) =
 -2 \mu;
\end{equation}
\begin{equation}\label{R5}
\lim_{k\to\infty} 
 \frac{E_a(k)}{\sigma(k)}= 
 {\mu_L+\mu}, \,\,
\lim_{k\to\infty}d(k) = 0,
\,\, 
\lim_{k\to\infty}c(k) =
-2 \mu .
\end{equation}
 
 \hfill $\bullet$ 
}

\vspace{0.25cm}

{\bf Proposition 2.}
{\it 
 The eigenfunctions 
(\ref{D1}),  with $A$ and $B$
given by the system 
(\ref{F3A}),  
do not check 
the compatibility condition (\ref{F-4}).

Proof. The relations 
(\ref{COMP}), (\ref{F-4}),
 (\ref{F-5})   give us
\begin{equation}\label{R9}
\sigma(k) = 
\frac{E_b(k)}
{\mu F(k,b)+ \mu_R} = \frac{E_a(k)}
{\mu_L- \mu F(k,a)},
\quad \forall k \geq 0.
\end{equation}
From
  (\ref{R4}), (\ref{R5}),
 (\ref{R9}) we obtain the existence of the limits 
$ \lim_{k\to\infty}
F(k,b)$,
$ \lim_{k\to\infty}
F(k,a)$.   We use the notations below and get
 $$  
 \lim_{k\to\infty}
F(k,b)=F(b),    \quad
\lim_{k\to\infty}
F(k,a)=F(a),            $$
\begin{equation}\label{R9A}
 F(b) \neq 1 
 \Rightarrow
\lim_{k\to\infty}
\frac{E_b(k)}{\sigma(k)}
\neq \mu+\mu_R;
\end{equation}
\begin{equation}\label{R9B}
 F(a) \neq - 1 
 \Rightarrow
\lim_{k\to\infty}
\frac{E_a(k)}{\sigma(k)}
\neq \mu_L+\mu.
\end{equation}
The last two relations
are in contradiction with
both possible relations (\ref{R4}),(\ref{R5}). 
Thus   we obtain
$F(b)=1, \quad F(b)=-1$.
The relation 
(\ref{R9}) for large $k$ 
gives us
\begin{equation}\label{T}
\frac{T(b)}{T(a)}=
\frac{\mu+\mu_R}
{\mu_L+\mu }.         
\end{equation}
Therefore we obtain 
the following unexpected restrictions,
which were not initially imposed: 
\begin{equation}\label{R10}
\frac{T(b)}{T(a)} =
\frac{\mu_R+\mu}{\mu_L+\mu} \Leftrightarrow
\mu=
\frac{\mu_RT(a)-\mu_LT(b)}{T(b)-T(a)},
\end{equation} 
\begin{equation}\label{R11}
0< \mu_L <
[\mu_RT(a)-\mu_LT(b)]/
[T(b)-T(a)]  < \mu_R.
\end{equation}                    
 The physical significance (related to fluid mechanics) of these  restrictions  
 is not clear. 
 Moreover, if we assume $T(b)< T(a)$, we need
$\mu_RT(a)< \mu_LT(b)$, 
so $\mu_R<\mu_L$.
Thus, in general,  the condition  
(\ref{F-4}) is not fulfilled.
}
\hfill $\bullet$

\vspace{0.5cm}

{\bf Remark 1.}
{\it
The 
eigenfunctions 
(\ref{D1})
 do not check 
the compatibility
condition (\ref{F-4}).
So the growth rates of the problem 
(\ref{J2A})-(\ref{COMP})
do not exist.
In fact,  in general,  
the problem
(\ref{J2A})-(\ref{COMP})
 doesn't make sense. 
However, in \cite{DA0} are given some estimations for 
the growth rates $\sigma$,
by using relation 
(\ref{F-3}).
Moreover, the results obtained in 
 \cite{DA0} are
  used in the papers
  \cite{DA1} - \cite{DA5}. 
We proved in \cite{PA} 
that the multi-layer model with constant
viscosities is not useful for minimizing the Saffman-Taylor instability, when 
the coefficients $A,B$ in  (\ref{D1}) are absolute constants.
The present paper 
can be considere  an improvement of  \cite{PA} 
for the case
$A=A(k), B=B(k)$.
}
\hfill $\bullet$

{\it 
\hspace{4cm} 
Simion Stoilow Institute  of Mathematics of the

\hspace{6cm} 
Romanian Academy

\hspace{4.4cm} 
 Calea Grivitei 21, Bucharest S1, Romania 
 }
 
 \hspace{5.6cm} 
 e-mail:   gelu.pasa@imar.ro

\end{document}